\documentclass[12pt]{amsart}
\usepackage[margin=30mm]{geometry}
\usepackage[all]{xy}
\usepackage{amssymb}
\usepackage{amsmath}
\usepackage{bigints}
\usepackage{bbm}
\usepackage{graphics,graphicx}
\usepackage{esint}
\usepackage{stfloats}
\usepackage{eepic}
\usepackage{epsfig}
\usepackage{breqn}
\usepackage{stmaryrd}
\usepackage{tikz}
\usepackage{hyperref}
\usepackage{wrapfig}
\usepackage{enumerate}

\usepackage{enumitem}

\newtheorem{theorem}{Theorem}[section]

\newtheorem{lemma}[theorem]{Lemma}

\newtheorem{remark}[theorem]{Remark}

\newcommand{\R}{\mathbb{R}}

\newcommand{\Om}{\Omega}
\newcommand{\om}{\omega}
\newcommand{\fii}{\varphi}
\newcommand{\eps}{\varepsilon}

\usepackage{caption} 
\usepackage{tikz}
\newcounter{mnotecount}[section]

\newcommand{\rmnote}[1]{}

\begin{document}
\title{Sensor Placement via Large Deviations in the Eikonal Equation}

\author[]{Ilias Ftouhi}

\address[Ilias Ftouhi]{Laboratoire de Math\'ematiques, Informatiques, Physique et Applications, Nîmes University, Site des Carmes, Place Gabriel P\'eri, 30000 Nîmes, France}
\email{ilias.ftouhi@unimes.fr}

\author[]{Enrique Zuazua}

\address[Enrique Zuazua]{\textsc{[1] Chair for Dynamics, Control, Machine Learning and Numerics, Alexander von Humboldt-Professorship, Department of Mathematics,  Friedrich-Alexander-Universit\"at Erlangen-N\"urnberg,91058 Erlangen, Germany.\vspace{1mm}}\newline
\textsc{[2] Chair of Computational Mathematics, University of Deusto, Av. de las Universidades, 24,
48007 Bilbao, Basque Country, Spain \vspace{1mm}}\newline
\textsc{[3] Departamento de Matem\'{a}ticas, Universidad Aut\'{o}noma de Madrid, 28049 Madrid, Spain.}\vspace{1mm}}
\email{enrique.zuazua@fau.de}

\date{\today}

\begin{abstract}
In this work, we address the problem of optimally placing a finite number of sensors within a given region so as to minimize the mean or maximal distance to the points of the domain. To tackle this natural geometric performance criterion, formulated in terms of distance functions, we combine tools from geometric analysis with a classical result of Varadhan~\cite{varadhan}, which provides an efficient approximation of the distance function via the solution of a simple elliptic PDE. The effectiveness of the proposed approach is demonstrated through illustrative numerical simulations.

\end{abstract}



\maketitle

\section{Introduction}
The optimal placement and design of sensors is a recurring challenge in industrial and applied contexts, ranging from urban planning to the monitoring of temperature and pressure in gas networks. In essence, sensors are optimally designed when they provide the most effective observation of the phenomenon under study. Typically, this design process is driven by precise objectives and shaped by constraints that are naturally formulated through partial differential equations, which capture the underlying physics of the system.

This article contributes to this line of research by presenting new insights into optimal design strategies in settings governed exclusively by geometric considerations. In particular, we investigate the problem of sensor placement in a purely geometric framework, without reference to a specific PDE model. The performance criteria we adopt are simple and natural, being defined in terms of distance functions. As we shall see, however, addressing this apparently elementary problem requires the use of analytical methods and naturally leads to the consideration of classical PDE-based approximations of distance functions.

These problems can then be formulated in a shape optimization framework. Indeed, given a set {$\Omega\subset\R^n$,} and a mass fraction $c\in (0,|\Om|)$, we state the following problem
$$\inf\{ \sup_{x\in\Om} d(x,\om)\ |\ |\om|=c\ \text{and}\ \om\subset \Om\},$$
where $d(x,\om):= \inf_{y\in \om}\|x-y\|$ is the minimal distance from $x$ to $\om$. 

In this way, our goal is to determine the optimal placement of sensors occupying a fraction 
$c \in (0,|\Omega|)$ of the available space in $\Omega$, where the physical process under consideration evolves, to minimize the maximal distance 
from any point in the domain to the nearest sensor.

This problem can be written in terms of the classical Hausdorff distance $d^H$, as  when $\om\subset \Om$, one has 
$$\sup_{x\in\Om} d(x,\om) = d^H(\om,\Om).$$

Thus, we consider the following problem 
\begin{equation}\label{problem:main_problem}
\inf \{d^H(\om,\Om) \ |\ |\om|=c\ \text{and}\ \om\subset \Om\},    
\end{equation}
where $c\in(0,|\Om|)$. 

By employing a homogenization strategy, in which the mass of many micro-sensors is uniformly 
distributed throughout $\Omega$, one can show that Problem~\eqref{problem:main_problem} does not 
admit an actual solution. Indeed, the infimum of the cost functional vanishes and is asymptotically 
approached by a sequence of disconnected sets with an increasing number of connected components 
densely covering $\Omega$. {Note that, in the absence of additional geometric constraints on the sensor set $\omega$, one can always construct sensors of zero measure whose Hausdorff distance to $\Omega$ is zero, i.e., $d^{H}(\omega,\Omega)=0$. For instance, taking $\omega=\Omega\cap\mathbb{Q}^d$ yields a set of Lebesgue measure zero that is dense in $\Omega$, and hence satisfies $d^{H}(\omega,\Omega)=0$. By considering sufficiently small neighborhoods around arbitrarily distributed rational points, this construction can be extended to sets of any prescribed measure. }

It is then mandatory to consider additional constraints on the sensor $\om$. For example, in \cite{FATTAH2025113866,f_zuazua}, the authors considered a convexity constraint on the sensor $\omega$, which ensures the existence of a minimizer.   

In the present paper, we consider the problem of minimizing the farthest distance to $N$ spherical sensors of a given radius $r$. The problem is formulated as follows 
\begin{equation}\label{prob:balls}
    \inf\{d^H(\cup_{i=1}^N B_i, \Om)\ |\ \forall i\in\{ 1,\dots,N \},\ \ \ B_i\in\Om\},
\end{equation}
where $B_i$ are balls of radius $r$. We then denote by $(x_i)$ the centers of the balls and write \eqref{prob:balls} as a finite dimensional optimization problem 
\begin{equation}\label{prob:f}
 \inf \{f(x_1,\dots,x_N)\ |\ x_i,\dots,x_N\in \mathcal{V}_{r}\},   
\end{equation}
where 
\begin{equation}\label{eq:v_r}
    \mathcal{V}_{r} := \{(x_1,\dots,x_N)\ |\ \forall i\ne j,\ \|x_i-x_j\|\ge 2r\ \text{and}\ \forall i,\ d(x_i,\partial \Omega)\ge r\}  
\end{equation}
and $$f(x_1,\dots,x_N) := \|d(\cdot,\cup_{i=1}^N B_i)\|_\infty 
 = d^H(\cup_{i=1}^N B_i, \Om).$$ 

 The set $\mathcal{V}_{r}$ is the set of centers $(x_i)$ such that the open sensors $B_i$ are all two by two disjoint and included in $\Omega$. 

This is the type of formulation one typically encounters in practical applications, since the shape and number of available sensors are often prescribed a priori.

Since the infinity norm $\|\cdot\|_\infty$ involved in the formulation of this problem is non-differentiable, it is generally impractical for numerical optimization, where accurate gradient formulas for the objective function are required. 

To avoid this technical difficulty, we work with the following regularized functional that involves the $p$-norm $\|\cdot\|_p$, with $p>1$. We are led to consider the minimization of  functionals of the form
 $$f_p(x_1,\dots,x_N):= \|d(\cdot,\cup_{i=1}^N B_i)\|_p = \left(\int_\Om d(x,\cup_{i=1}^N B_i)^p dx\right)^{1/p},$$
 known as \textit{average distance} functionals. We refer to \cite{zbMATH05349205,zbMATH06110561} for reviews of average distance problems and relevant references. These problems are also known in the operation research community under the name of \textit{the continuous Fermat--Weber Problems}, see for example \cite{zbMATH05580335} and the references therein. 

{The main contribution of this paper is the use of a classical result of Varadhan~\cite{varadhan}, stated in Theorem \ref{th:varadhan}, to approximate the distance function through the solution of a suitable elliptic PDE. This approach naturally leads to a family of approximated problems (see Section~\ref{s:varadhan} and more precisely Lemma \ref{lem:approx}), for which we perform a sensitivity analysis using shape-derivative techniques, see Section~\ref{s:derivative} where a shape derivation result is presented in Theorem \ref{prop:shape_derivative}. Then, these tools are applied efficiently to numerical optimization, producing the numerical results presented in Section~\ref{s:numerics}.}

\section{State of the art and Applications}
In this section, we present a non-exhaustive overview of the state of the art concerning the optimal placement of sensors.

A first example arises in gas networks, where strategically positioning sensors is crucial for the efficient distribution and management of resources such as natural gas and hydrogen. 
Sensors are typically deployed at critical locations, such as pipelines, compressor stations, and distribution hubs to monitor key parameters including pressure, temperature, flow rates, and gas composition. Such deployments enable the rapid detection of leaks or pressure fluctuations, allowing operators to respond promptly and thereby mitigate safety risks and environmental damage. Moreover, the collected data supports performance optimization, reduces energy losses, and enhances the overall reliability and cost-effectiveness of the network.

Another illustrative application is found in urban planning, where the optimal placement of facilities within cities plays a decisive role in improving living conditions. Careful design of public spaces, transportation networks, healthcare facilities, and green areas improves access to essential services, alleviates congestion, and promotes sustainability. By accounting for factors such as population density, accessibility, and environmental impact, urban planners can foster more livable, efficient, and resilient cities, thereby improving overall quality of life.

The problem of optimal sensor (or facility) placement has thus been extensively investigated across a wide range of domains, underlining both its practical relevance and theoretical significance. Below, we provide a representative selection of related problems and references.

\begin{itemize}
      \item The problem of designing actuators within a prescribed set, to optimize the control performance of a process, has been extensively investigated in the framework of optimal control. Several methodologies for effective actuator placement have been developed, notably those based on spectral techniques, as discussed in~\cite{PrivatTrelatZuazua13,zbMATH06440180,zbMATH06722979}, with applications primarily to control problems for wave and heat processes. {We also refer to \cite{zbMATH08136281,zbMATH07745053,zbMATH07984160,zbMATH07898995} for a non-exhaustive list of references linking theory with applications.}
      
  \item The problem of minimizing the average distance within a given subset has long been a subject of interest in the optimal transport community~\cite{zbMATH01996419,zbMATH05058778}. Comprehensive accounts of its development, together with surveys of the field, can be found in~\cite{zbMATH05349205,zbMATH06110561}. More recently, attention has shifted toward minimizing the maximal distance, thereby broadening both the theoretical foundations and the range of practical applications of the problem~\cite{zbMATH,zbMATH06996635,f_zuazua,zbMATH06919706}.

\item The problem of optimally positioning cavities within a domain to control specific eigenvalues of a differential operator has been the subject of extensive research. A variety of approaches have been proposed to minimize or maximize such eigenvalues, giving rise to a rich and diverse literature. For a comprehensive and up-to-date review, including a detailed introduction and an extensive bibliography, we refer the reader to~\cite{ftouhi_steklov}.

\item This perspective has been applied to mathematical biology problems, such as optimizing resource distribution to maximize population size \cite{zbMATH07155098,zbMATH07168708}, as well as the mathematical analysis of optimal habitat configurations for species persistence \cite{zbMATH05240724}. Also, optimal sensor placement plays a crucial role in the study of the brain, as it directly determines the quality, interpretability, and efficiency of the data acquired from neuroimaging and electrophysiological techniques such as electroencephalogram (EEG) and near-infrared spectroscopy (NIRS), see for example \cite{Choi2006EEGReference,Rogers2025Colocalized} . 
\end{itemize}

The wide range of applications and the theoretical complexity of these problems highlight the need for developing robust mathematical methods for the optimal design and placement of actuators.

It is important to note, however, that in this paper we approach the problem from a purely geometric perspective, without reference to any specific application. In this sense, the sensor location problems we address have a form of universal validity, as they rely solely on geometric considerations. The outcomes of our methodologies can thus be viewed as an initialization step, to be further refined when incorporating the specific knowledge of a given application--for instance, by constraining the admissible functions to be solutions of a particular partial differential equation, as is often the case in practice.

\section{Approximating the distance function via a Varadhan's result}\label{s:varadhan}

In the present paper, we are considering simple and natural criteria of performance based on distance functions. More precisely, we are interested in finding the optimal placement of sensors inside a given region in such a way as to minimize the maximum and mean distance to the points of the domain. The prevailing approach to distance computation is to solve the
eikonal equation
\begin{equation}\label{prob:eikonal}
   |\nabla d| = 1 
\end{equation}
with Dirichlet boundary conditions $d = 0$ on the boundary. 

The Eikonal equation is known to present several well-documented analytical and computational challenges:
\begin{itemize}
\item Nonlinearity and non-smoothness: The equation is highly nonlinear and its solutions are typically only Lipschitz continuous, with corners or ridges where classical derivatives do not exist.
\item 	Singularities and shocks: The solution may develop singularities (e.g., at points equidistant from multiple boundaries, the cut locus of the boundary), where $\nabla u$ is not well-defined.
	\item	Viscosity solutions: Because classical solutions often fail to exist, one must work in the framework of viscosity solutions, which requires special numerical schemes.

	\item 	Numerical challenges: Standard discretizations can produce instabilities or fail to converge. Specialized methods like the fast marching method and the fast sweeping method are commonly used to handle the geometric structure efficiently.
	\item	Boundary sensitivity: The solution strongly depends on boundary data, and inaccuracies in its discretization can propagate throughout the domain.
\end{itemize}	
	
To circumvent these technical difficulties, and following the ideas in~\cite{crane}, we propose to approximate distance functions using a classical PDE result due to Varadhan~\cite{varadhan}:

\begin{theorem}\label{th:varadhan}
Let $\Omega$ be an open subset of $\R^n$ and $\varepsilon>0$, we consider the problem 
\begin{equation}\label{prob:varadhan}
    \begin{cases}
    w_\eps-\eps \Delta w_\eps=0 &\qquad\mbox{in $\Omega,$}\vspace{1mm} \\
   \ \ \ \ \ \ \ \ \ \  w_\eps = 1& \qquad\mbox{on $\partial \Omega$.} 
    \end{cases}
\end{equation}

We have
$$\lim_{\eps\rightarrow 0} -\sqrt{\eps}\ln{w_\eps(x)} = d(x,\partial \Omega):= \inf_{y\in \partial \Omega} \|x-y\|,$$
uniformly over compact subsets of $\Om$.  
\end{theorem}

The original proof of Theorem~\ref{th:varadhan} relies on the maximum principle together with explicit solutions of equation~\eqref{prob:varadhan} in the case of balls. This approach, however, does not provide a convergence rate. Since the problem~\eqref{prob:varadhan} is closely connected to the distance function, it is natural to expect a relation with the Eikonal equation~\eqref{prob:eikonal}. This connection was fully clarified following the development of the theory of viscosity solutions, introduced in the early 1980s by Pierre-Louis Lions and Michael G. Crandall~\cite{zbMATH03966825}. In fact, one can readily verify that the function $v_\eps:= -\sqrt{\eps}\log{w_\eps}$ satisfies the following viscous approximation of the Eikonal equation:
\begin{equation}\label{prob:viscous_eikonal}
    \begin{cases}
    1-|\nabla v_\eps|^2 +\sqrt{\eps}\Delta v_\eps = 0 &\qquad\mbox{in $\Omega,$}\vspace{1mm} \\
   \ \ \ \ \ \ \ \ \ \ \ \ \ \ \ \ \ \ \ \ \ \  v_\eps = 0& \qquad\mbox{on $\partial \Omega$.} 
    \end{cases}
\end{equation}

Such an equation has been considered in \cite{japonais}, where, by using a doubling variables method,  the authors prove in Theorem 3.1   that there exists a constant $C(\Omega)$ depending only on $\Omega$ such that 
\begin{equation}\label{eq:rate}
    \forall x\in\Omega, \forall \eps > 0,\ \ \ |v_\eps(x)-d(x,\partial \Omega)|\leq C(\Omega)\eps^{\frac{1}{4}}.
\end{equation}

We now show how this result can be employed to construct an approximation 
of the objective function $f$ defined in problem~\eqref{prob:f}. 
The key idea is to embed $\Omega$ into a sufficiently large box $D$, 
chosen so that, for every point in $\Omega$, the distance  to the exterior boundary $\partial D$ is smaller than the distance to the sensors 
$\cup_{i=1}^N B_i$.

\begin{figure}[!h]
    \centering

\tikzset{every picture/.style={line width=0.75pt}} 

\begin{tikzpicture}[x=0.75pt,y=0.75pt,yscale=-1,xscale=1]

\draw  [line width=2.25]  (218.33,91) .. controls (242.76,75.26) and (410.33,109) .. (412.33,135) .. controls (414.33,161) and (315.33,191) .. (376.33,239) .. controls (437.33,287) and (163,296) .. (143,266) .. controls (123,236) and (170.91,193.74) .. (177.91,170.74) .. controls (184.91,147.74) and (193.91,106.74) .. (218.33,91) -- cycle ;
\draw  [fill={rgb, 255:red, 80; green, 227; blue, 194 }  ,fill opacity=1 ] (210.16,120.14) .. controls (210.17,114.56) and (214.71,110.06) .. (220.28,110.08) .. controls (225.85,110.1) and (230.35,114.63) .. (230.34,120.2) .. controls (230.32,125.77) and (225.79,130.28) .. (220.21,130.26) .. controls (214.64,130.24) and (210.14,125.71) .. (210.16,120.14) -- cycle ;
\draw  [fill={rgb, 255:red, 80; green, 227; blue, 194 }  ,fill opacity=1 ] (200,184.17) .. controls (200,178.55) and (204.55,174) .. (210.17,174) .. controls (215.78,174) and (220.33,178.55) .. (220.33,184.17) .. controls (220.33,189.78) and (215.78,194.33) .. (210.17,194.33) .. controls (204.55,194.33) and (200,189.78) .. (200,184.17) -- cycle ;
\draw  [fill={rgb, 255:red, 80; green, 227; blue, 194 }  ,fill opacity=1 ] (206.33,261.17) .. controls (206.33,255.92) and (210.59,251.67) .. (215.83,251.67) .. controls (221.08,251.67) and (225.33,255.92) .. (225.33,261.17) .. controls (225.33,266.41) and (221.08,270.67) .. (215.83,270.67) .. controls (210.59,270.67) and (206.33,266.41) .. (206.33,261.17) -- cycle ;
\draw  [fill={rgb, 255:red, 80; green, 227; blue, 194 }  ,fill opacity=1 ] (250,189.67) .. controls (250,184.33) and (254.33,180) .. (259.67,180) .. controls (265.01,180) and (269.33,184.33) .. (269.33,189.67) .. controls (269.33,195.01) and (265.01,199.33) .. (259.67,199.33) .. controls (254.33,199.33) and (250,195.01) .. (250,189.67) -- cycle ;
\draw  [fill={rgb, 255:red, 80; green, 227; blue, 194 }  ,fill opacity=1 ] (358,137.67) .. controls (358,132.33) and (362.33,128) .. (367.67,128) .. controls (373.01,128) and (377.33,132.33) .. (377.33,137.67) .. controls (377.33,143.01) and (373.01,147.33) .. (367.67,147.33) .. controls (362.33,147.33) and (358,143.01) .. (358,137.67) -- cycle ;
\draw  [fill={rgb, 255:red, 80; green, 227; blue, 194 }  ,fill opacity=1 ] (300,250.17) .. controls (300,244.55) and (304.55,240) .. (310.17,240) .. controls (315.78,240) and (320.33,244.55) .. (320.33,250.17) .. controls (320.33,255.78) and (315.78,260.33) .. (310.17,260.33) .. controls (304.55,260.33) and (300,255.78) .. (300,250.17) -- cycle ;
\draw  [fill={rgb, 255:red, 80; green, 227; blue, 194 }  ,fill opacity=1 ] (313,177.17) .. controls (313,171.55) and (317.55,167) .. (323.17,167) .. controls (328.78,167) and (333.33,171.55) .. (333.33,177.17) .. controls (333.33,182.78) and (328.78,187.33) .. (323.17,187.33) .. controls (317.55,187.33) and (313,182.78) .. (313,177.17) -- cycle ;
\draw  [fill={rgb, 255:red, 80; green, 227; blue, 194 }  ,fill opacity=1 ] (280,110.67) .. controls (280,105.33) and (284.33,101) .. (289.67,101) .. controls (295.01,101) and (299.33,105.33) .. (299.33,110.67) .. controls (299.33,116.01) and (295.01,120.33) .. (289.67,120.33) .. controls (284.33,120.33) and (280,116.01) .. (280,110.67) -- cycle ;
\draw  [color={rgb, 255:red, 208; green, 2; blue, 27 }  ,draw opacity=1 ][dash pattern={on 6.75pt off 4.5pt}][line width=2.25]  (57,57.33) -- (494.33,57.33) -- (494.33,321.33) -- (57,321.33) -- cycle ;
\draw [line width=1.5]  [dash pattern={on 1.69pt off 2.76pt}]  (177.91,170.74) -- (200.91,178.74) ;
\draw [shift={(177.91,170.74)}, rotate = 19.18] [color={rgb, 255:red, 0; green, 0; blue, 0 }  ][fill={rgb, 255:red, 0; green, 0; blue, 0 }  ][line width=1.5]      (0, 0) circle [x radius= 4.36, y radius= 4.36]   ;
\draw    (188.41,181.74) -- (184.34,205.03) ;
\draw [shift={(184,207)}, rotate = 279.9] [color={rgb, 255:red, 0; green, 0; blue, 0 }  ][line width=0.75]    (10.93,-3.29) .. controls (6.95,-1.4) and (3.31,-0.3) .. (0,0) .. controls (3.31,0.3) and (6.95,1.4) .. (10.93,3.29)   ;
\draw [color={rgb, 255:red, 74; green, 144; blue, 226 }  ,draw opacity=1 ][line width=1.5]  [dash pattern={on 1.69pt off 2.76pt}]  (58,170) -- (177.91,170.74) ;
\draw [color={rgb, 255:red, 74; green, 144; blue, 226 }  ,draw opacity=1 ]   (110,124) -- (115.7,161.02) ;
\draw [shift={(116,163)}, rotate = 261.25] [color={rgb, 255:red, 74; green, 144; blue, 226 }  ,draw opacity=1 ][line width=0.75]    (10.93,-3.29) .. controls (6.95,-1.4) and (3.31,-0.3) .. (0,0) .. controls (3.31,0.3) and (6.95,1.4) .. (10.93,3.29)   ;

\draw (385,189.4) node [anchor=north west][inner sep=0.75pt]  [font=\Large]  {$\Omega $};
\draw (502,121.4) node [anchor=north west][inner sep=0.75pt]  [font=\LARGE]  {$\textcolor[rgb]{0.82,0.01,0.11}{D}$};
\draw (233,109.4) node [anchor=north west][inner sep=0.75pt]  [font=\small]  {$\textcolor[rgb]{0.25,0.46,0.02}{B}\textcolor[rgb]{0.25,0.46,0.02}{_{1}}$};
\draw (302,100.4) node [anchor=north west][inner sep=0.75pt]  [font=\small]  {$\textcolor[rgb]{0.25,0.46,0.02}{B}\textcolor[rgb]{0.25,0.46,0.02}{_{2}}$};
\draw (381,128.4) node [anchor=north west][inner sep=0.75pt]  [font=\small]  {$\textcolor[rgb]{0.25,0.46,0.02}{B}\textcolor[rgb]{0.25,0.46,0.02}{_{3}}$};
\draw (221,165.4) node [anchor=north west][inner sep=0.75pt]  [font=\small]  {$\textcolor[rgb]{0.25,0.46,0.02}{B}\textcolor[rgb]{0.25,0.46,0.02}{_{4}}$};
\draw (273,181.4) node [anchor=north west][inner sep=0.75pt]  [font=\small]  {$\textcolor[rgb]{0.25,0.46,0.02}{B}\textcolor[rgb]{0.25,0.46,0.02}{_{5}}$};
\draw (335,168.4) node [anchor=north west][inner sep=0.75pt]  [font=\small]  {$\textcolor[rgb]{0.25,0.46,0.02}{B}\textcolor[rgb]{0.25,0.46,0.02}{_{6}}$};
\draw (224,239.4) node [anchor=north west][inner sep=0.75pt]  [font=\small]  {$\textcolor[rgb]{0.25,0.46,0.02}{B}\textcolor[rgb]{0.25,0.46,0.02}{_{7}}$};
\draw (322,240.4) node [anchor=north west][inner sep=0.75pt]  [font=\small]  {$\textcolor[rgb]{0.25,0.46,0.02}{B}\textcolor[rgb]{0.25,0.46,0.02}{_{8}}$};
\draw (163,210.4) node [anchor=north west][inner sep=0.75pt]  [font=\small]  {$d( B_{4} ,\partial \Omega )$};
\draw (76,103.4) node [anchor=north west][inner sep=0.75pt]  [font=\small,color={rgb, 255:red, 74; green, 144; blue, 226 }  ,opacity=1 ]  {$d( M_{4} ,\partial D)$};
\draw (156,147.4) node [anchor=north west][inner sep=0.75pt]  [font=\small]  {$M_{4}$};

\end{tikzpicture}

    \caption{The box $D$ containing the domain $\Omega$ and the sensors $(B_i)$. As it can be seen in the figure, the box $D$ is chosen such that $d(B_4,\partial \Omega)=d(B_4,M_4)\leq d(M_4,\partial D)$}
    \label{fig:box}
\end{figure}
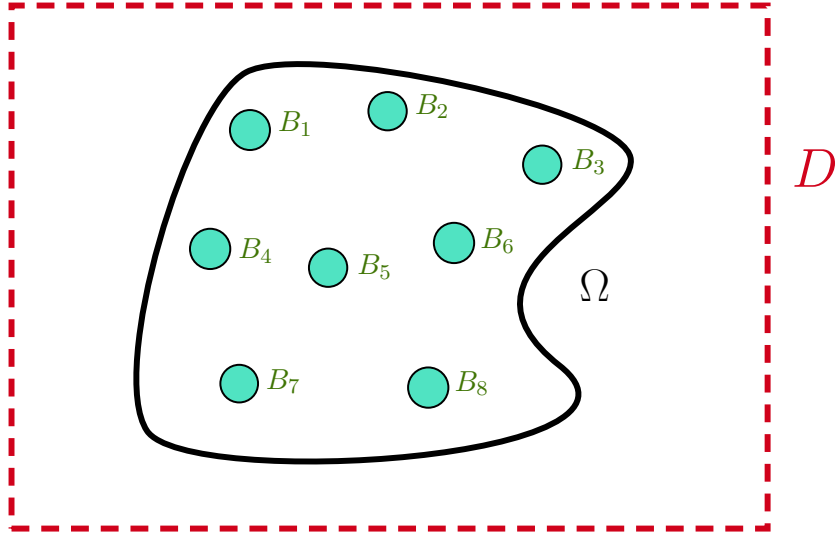

\begin{lemma}\label{lem:approx}
Let $\Omega$ be an open bounded domain of $\R^n$. We consider $N$ disjoint sensors $B_i\subset \Omega$ and $D$ a sufficiently large box containing $\Omega$ such that for every point of $\Omega$ the distance to the sensors $\cup_{i=1}^N B_i$ is smaller than the distance to the boundary of the box $\partial D$, see Figure \ref{fig:box}.  

For every $\eps>0$, we consider the problem
\begin{equation*}
    \begin{cases}
    w_\eps-\eps \Delta w_\eps=0 &\qquad\mbox{in $D\backslash\cup_{i=1}^N B_i$,}\vspace{1mm} \\
   \ \ \ \ \ \ \ \ \ \  w_\eps = 1& \qquad\mbox{on $ \cup_{i=1}^N \partial B_i\cup \partial D$.} 
    \end{cases}
\end{equation*}

The function $v_\eps:x\longmapsto -\sqrt{\eps}\ln{w_\eps(x)}$ uniformly converges to $$d(\cdot,\cup_{i=1}^N B_i): x\longmapsto \inf\limits_{y\in \cup_{i=1}^N B_i} \|x-y\|$$ on $\Om\backslash \cup_{i=1}^N B_i$. More precisely, there exists a constant $C$ depending only on $D$ and $(B_i)$ such that 
\begin{equation*}
    \forall \eps > 0,\forall x\in\Omega\backslash \cup_{i=1}^N B_i, \ \ \ |v_\eps(x)-d(\cdot,\cup_{i=1}^N B_i)|\leq C\eps^{\frac{1}{4}}.
\end{equation*}
\end{lemma}
\begin{proof}
By Theorem \ref{th:varadhan}, the function $v_\eps$ uniformly converges on $D\backslash\cup_{i=1}^N B_i$ to the distance function 
$$d: x\longmapsto \inf \{\|x-y\|\ |\ y\in  \cup_{i=1}^N B_i \cup \partial D\}.$$  

By the choice of a sufficiently large box $D$, we have for every $x\in \Om\backslash \cup_{i=1}^N B_i$
$$ d(x,\cup_{i=1}^N \partial B_i\cup \partial D) =  \min(d(x,\cup_{i=1}^N B_i),d(x,\partial D))= d(x,\cup_{i=1}^N B_i).$$

The quantitative estimate is a direct consequence of \eqref{eq:rate}. 
\end{proof}

Now that we have an approximation $v_\eps$ of the distance to the sensors $B_i$ via Varadhan's result, we consider the following approximation of $f$: 
$$f_{p,\eps}: (x_1,\dots,x_N)\in \mathcal{V}_{r}\longmapsto \|v_\eps\|_p =  \left(\int_\Om v_\eps^p dx\right)^{1/p},$$
where $\mathcal{V}_{r}$ is defined in \eqref{eq:v_r}. And, consequently,  we focus on the numerical resolution of problems
\begin{equation}\label{prob:f_approx}
 \inf \{f_{p,\eps}(x_1,\dots,x_N)\ |\ (x_1,\dots,x_N)\in \mathcal{V}_{r}\},   
\end{equation}
with $p>0$, or equivalently of the problem 
\begin{equation}\label{prob:f_approx_eq}
 \inf \{\int_{\Om\backslash \cup_{i=1}^N B_i} v_\eps^p dx\ |\ (x_1,\dots,x_N)\in \mathcal{V}_{r}\}.  
\end{equation}
The first step in determining the optimal configuration is to compute 
the gradient of the function $f_{p,\varepsilon}$. This is achieved by 
applying shape differentiation with respect to perturbations corresponding 
to translations of the sensors $B_i$. To derive efficient formulas for these shape derivatives, a suitable adjoint state is introduced. Once the 
gradients are available, they are employed in a gradient descent scheme 
to minimize $f_{p,\varepsilon}$ and thereby obtain the optimal placement 
of the sensors $B_i$.

{As we will see in remark \ref{rk:local}, although efficient and straightforward to implement, gradient methods suffer from a fundamental limitation since they only guarantee convergence to a critical point or a local minimum, when a second order optimality condition is used, not necessarily to the global one. 
Developing efficient methods that guarantee convergence to the global minimizer is a challenging problem and lies beyond the scope of this paper.}

\section{Numerical simulations}\label{s:numerics}
The present section is divided into two parts: first, we compute the shape derivative of the shape functionals approximating the average distance; then, we provide some numerical results on the optimal placement of sensors.

\subsection{Computation of the gradients}\label{s:derivative}
We recall that we are interested in minimizing the functions 
$$f_{p,\eps}(x_1,\dots,x_N):=  \left(\int_\Om (-\sqrt{\eps}\log w_\eps)^p dx\right)^{1/p},$$
or equivalently, minimizing the functions 
$$g_{p,\eps}(x_1,\dots,x_N):= f_{p,\eps}(x_1,\dots,x_N)^p =  \int_\Om (-\sqrt{\eps}\log w_\eps)^p dx,$$
with $p\ge 1$ and $\eps>0$ and $w_\eps$ is the solution of the problem  
\begin{equation*}
    \begin{cases}
    w_\eps-\eps \Delta w_\eps=0 &\qquad\mbox{in $D\backslash\cup_{i=1}^N B_i$,}\vspace{1mm} \\
   \ \ \ \ \ \ \ \ \ \ \  w_\eps = 1& \qquad\mbox{on $\partial \cup_{i=1}^N B_i\cup \partial D$,} 
    \end{cases}
\end{equation*}
where $D$ is a sufficiently large box containing $\Omega$ (see Lemma \ref{lem:approx}) and $(B_i)_i$ are disjoint balls of radius $r$ included in $\Omega$.  

To apply a gradient descent method, we have to efficiently compute the gradients of the functions $g_{p,\eps}$. To do so, we will use shape derivation techniques, where the perturbation of the centers of the sensors $B_i$ will be seen as a shape perturbation of the boundary of the perforated domain $\Omega\backslash \cup_{i=1}^N B_i$. Our result in this direction is then stated as follows: 
\begin{theorem}\label{prop:shape_derivative}
    Let $p\ge 1$ and $\eps>0$. The function $g_{p,\eps}$ is differentiable on the interior of $\mathcal{V}_r$ (defined in \eqref{eq:v_r}) and we have for every $(i,k)\in \llbracket 1,N\rrbracket\times \llbracket 1,n\rrbracket$,
    \begin{equation}\label{eq:formula_particuliere}
        \frac{\partial g_{p,\eps}}{\partial x_i^k} = \int_{\partial B_i} \left((-\sqrt{\eps}\log(w_\eps))^p + \eps\frac{\partial w_k}{\partial \nu}\frac{\partial q_\eps}{\partial \nu} \right)\nu_k d\sigma,
    \end{equation}
    where $\nu = (\nu_1,\dots,\nu_n)\in \R^n$ is the normal vector to the boundary of $B_i$ and $q_\eps$ is the solution of the following adjoint problem
   $$ \begin{cases}
    q_\eps-\eps \Delta q_\eps=-\frac{p}{w_\eps}(-\sqrt{\eps}\log(w_\eps))^{p-1}\mathbbm{1}_\Om &\qquad\mbox{in $D\backslash\cup_{i=1}^N B_i$,}\vspace{1mm} \\
   \ \ \ \ \ \ \ \ \  q_\eps = 0 & \qquad\mbox{on $\cup_{i=1}^N\partial B_i\cup \partial D$.} 
    \end{cases} $$
\end{theorem}

\begin{proof}   
Let us take a general setting and consider the functional 
$$J(\omega):= \int_{D\backslash\om } j(u)\mathbbm{1}_\Om dx =  \int_{\Om\backslash\om } j(u)dx,$$ where $\om$ is a smooth subdomain of $\Omega$ and $j:\R\longrightarrow \R$ is a given function and $u$ is the solution of the problem 
\begin{equation*}
    \begin{cases}
    u-\eps \Delta u=-1 &\qquad\mbox{in $D\backslash\om$,}\vspace{1mm} \\
   \ \ \ \ \ \ \ \ \ \ u = 0& \qquad\mbox{on $\partial \om\cup \partial D$,} 
    \end{cases}
\end{equation*}
with $\eps>0$.

Consider $V:\R^n\rightarrow\R^n$ a smooth vector field supported in a small neighborhood of $\om$. Let us prove that the functional $J$ admits a first-order shape derivative in the direction $V$, and we have 
\begin{equation}\label{eq:derivative_general}
    J'(\om,V) = \int_{\partial \om} \left(j(u) + \eps\frac{\partial u}{\partial \nu}\frac{\partial q}{\partial \nu} \right)\langle V,\nu \rangle d\sigma,
\end{equation}
    where $q$ is the solution of the following adjoint problem 
$$
\begin{cases}
    q-\eps \Delta q=-j'(u)\mathbbm{1}_\Om &\qquad\mbox{in $D\backslash\om$,}\vspace{1mm} \\
   \ \ \ \ \ \ \ \ \ q = 0 & \qquad\mbox{on $\partial \om\cup \partial D$.} 
\end{cases}
$$

The shape derivative of $u$ with respect to the perturbation $V$ is denoted by $u'$ and satisfies the following equation   

$$
    \begin{cases}
    u'-\eps \Delta u'=-0 &\qquad\mbox{in $D\backslash\om$,}\vspace{1mm} \\
   \ \ \ \ \ \ \ \ \ \ u' = -\frac{\partial u}{\partial \nu}\langle V,\nu \rangle  & \qquad\mbox{on $\partial \om\cup \partial D$.} 
    \end{cases}    
$$
{The last claims are standard, we refer to \cite[Section 5.3]{HPb} for the formal proofs.} On the other hand, we have
$$J'(\om,V) = \int_{\partial (\Om\backslash\om)} j(u)\langle V,\nu \rangle d\sigma + \int_{\Om\backslash \om} j'(u)u' dx,$$
see for example \cite[Theorem 5.2.2]{HPb}. 

Since $V$ is supported on a small neighborhood of $\omega$, it vanishes on $\partial \Omega$. Therefore, the last equality can be written as follows
\begin{equation}\label{eq:first}
J'(\om,v) = \int_{\partial \om} j(u)\langle V,\nu \rangle d\sigma + \int_{\Om\backslash \om} j'(u)u' dx.
\end{equation}
We now show how the terms involving $u'$ can be eliminated. We multiply both sides of the adjoint equation by a function $\varphi\in H^1(D\backslash \om)$, then integrate on $D\backslash \om$. We have 
{$$\int_{D\backslash \om} q\varphi dx-\eps \int_{D\backslash \om}\Delta q \varphi dx = -\int_{D\backslash \om} j'(u)\varphi\mathbbm{1}_\Om dx = -\int_{\Om\backslash \om} j'(u)\varphi dx.$$}
By integration by parts, we have 
$$\int_{D\backslash \om}\Delta q \varphi dx = -\int_{D\backslash \om}\nabla q \nabla \varphi dx + \int_{\partial(D\backslash \om)}\frac{\partial q}{\partial \nu} \varphi d\sigma.$$
Therefore, we have 
$$\forall \fii\in  H^1(D\backslash \om),\ \ \ \int_{\Om\backslash \om} j'(u)\fii dx = -\int_{D\backslash \om} q\fii dx -\eps\int_{D\backslash \om} \nabla q \nabla u dx + \eps \int_{\partial D\cup \partial \om} \frac{\partial q}{\partial \nu} \fii d\sigma. $$
By taking $\fii = u'$ in the last equality, we have 
\begin{equation}\label{eq:second}
\int_{\Om\backslash \om} j'(u)u' dx = -\int_{D\backslash \om} q u' dx -\eps\int_{D\backslash \om} \nabla q \nabla u dx + \eps \int_{\partial D\cup \partial \om} \frac{\partial q}{\partial \nu} u' d\sigma.    
\end{equation}
Now, by considering the equation satisfied by $u'$, we have 
$$\forall \fii\in  H^1(D\backslash \om),\ \ \  \int_{D\backslash \om} u'\fii dx +\eps\int_{D\backslash \om} \nabla u' \nabla \fii dx - \eps \int_{\partial D\cup \partial \om} \frac{\partial u'}{\partial \nu} \fii d\sigma = 0. $$
Taking $\fii = q$, yields 
\begin{equation}\label{eq:third}
  \int_{D\backslash \om} u'q dx +\eps\int_{D\backslash \om} \nabla u' \nabla q dx = 0.  
\end{equation}

The claim of \eqref{eq:derivative_general} then follows from \eqref{eq:first}, \eqref{eq:second} and \eqref{eq:third}.

At last, Formula \eqref{eq:formula_particuliere} is obtained by using \eqref{eq:derivative_general} with $\omega = \cup_{i=1}^N B_i$, $V:x\longmapsto e_k$ for $k\in \llbracket 1,n \rrbracket$, with $(e_k)$ being the canonical basis of $\R^n$, and $j: x\longmapsto (\sqrt{\eps}\log{(x+1)})^p$. 
\end{proof}

\subsection{Numerical experiments}
In this section, we present a series of numerical simulations. The computations were carried out using a gradient descent algorithm, where the gradients were obtained from the result of Theorem~\ref{prop:shape_derivative}.

Gradient descent is a widely used optimization method that iteratively updates parameters in the direction of the negative gradient of a cost function. Although efficient and straightforward to implement, it suffers from a fundamental limitation: it guarantees convergence only to a local minimum, not necessarily to the global one. This issue is particularly acute for non-convex functions, which may exhibit multiple local minima and saddle points. As a purely local method, the outcome of gradient descent strongly depends on the initial starting point, and the algorithm may easily become trapped in suboptimal regions of the landscape.

To mitigate this difficulty, we adopted a multi-start strategy, running gradient descent from a variety of initial conditions. This increases the likelihood of exploring different regions of the search space and improves the chances of identifying better -- possibly global -- minima by comparing the results of multiple independent runs.

{The simulations of both methods were conducted using \texttt{Matlab} software on a personal \texttt{DELL} computer with a quadcore Intel Core i7-2600 with 16GB RAM. The resolution of PDE \eqref{prob:varadhan} was done by using finite element method via the toolbox \texttt{PDEtool} and with the parameter $\eps$ taken equal to $10^{-3}$. The optimization procedure was performed through \texttt{'fmincon'} routine with the following settings: 
\begin{itemize}
    \item \texttt{Algorithm:} \texttt{"sqp"}.
    \item \texttt{ConstraintTolerance: 1.0000e-06}.
    \item \texttt{OptimalityTolerance: 1.0000e-06}.
    \item \texttt{StepTolerance: 1.0000e-06}.
    \item \texttt{SpecifyObjectiveGradient: "on"}.
\end{itemize}}

\begin{figure}[ht]
    \centering
\includegraphics[scale=.6]{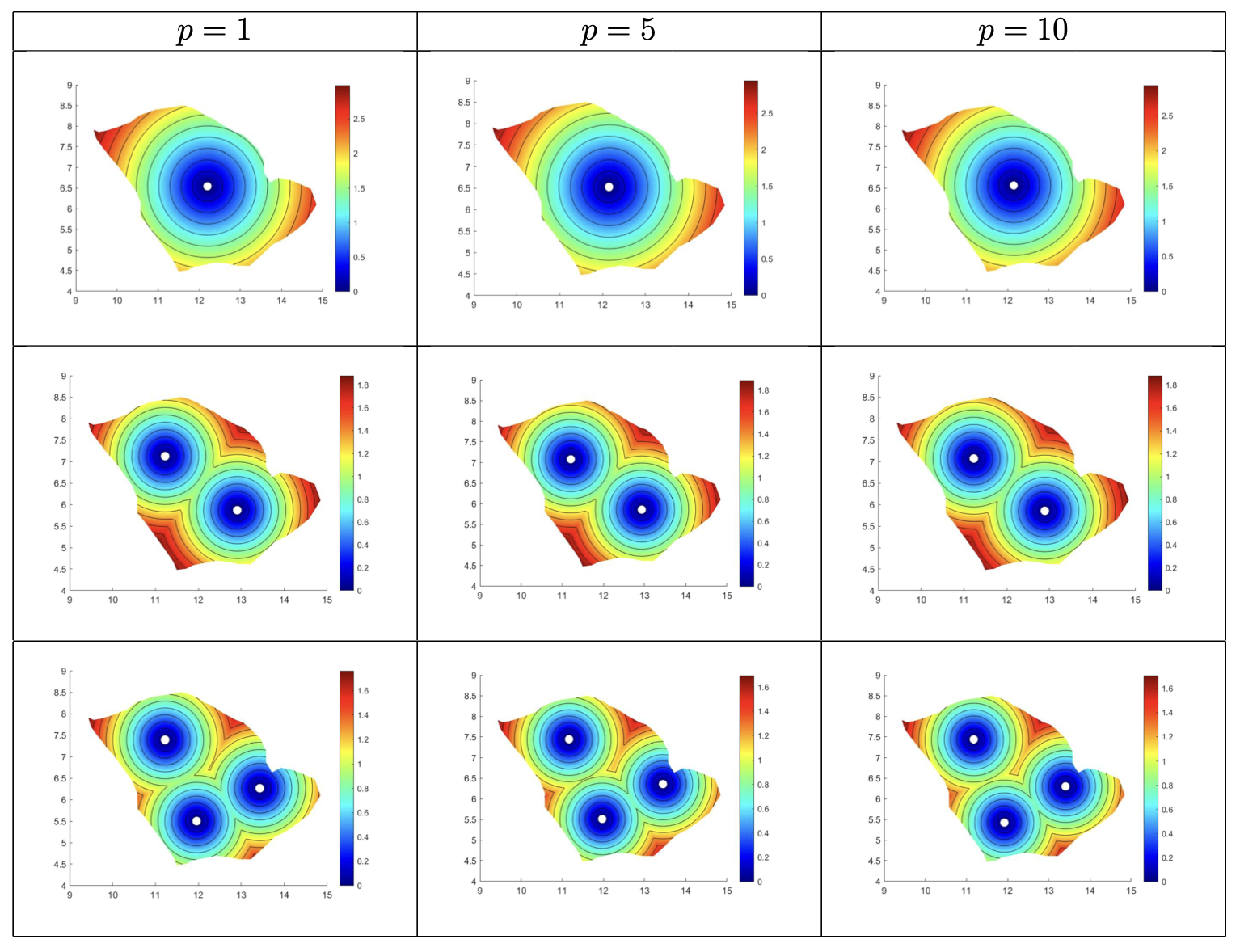}
    \caption{Optimal placement of $N\in \{1,2,3\}$ sensors inside the country of Saudi Arabia for different values of $p$.}
    \label{fig:sensors_algeria}
\end{figure}

\begin{figure}[ht]
    \centering
\includegraphics[scale=.7]{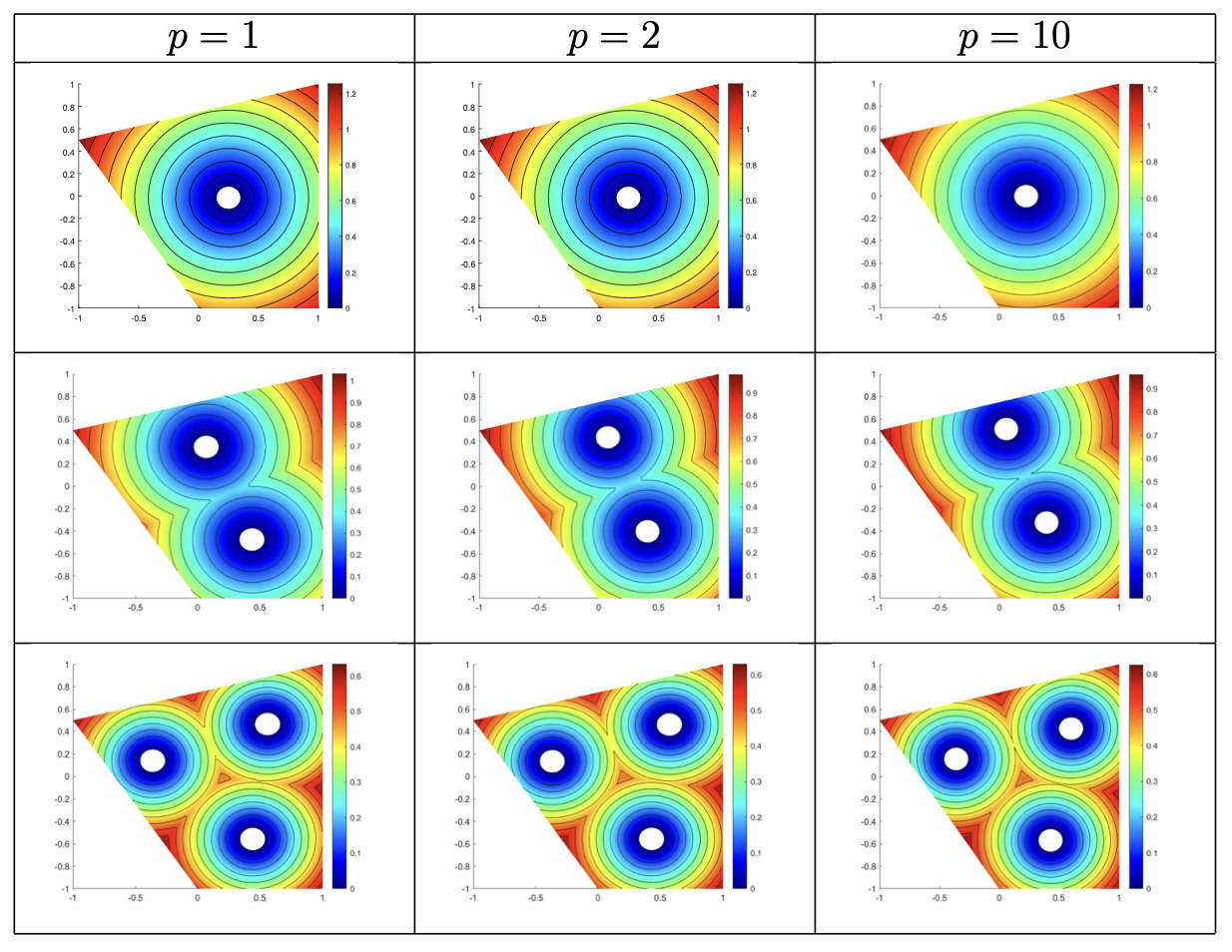}
    \caption{Optimal placement of $N\in \{1,2,3\}$ sensors in rhombus for different values of $p$.}
    \label{fig:sensors_3}
\end{figure}

\begin{figure}[ht]
    \centering
\includegraphics[scale=.55]{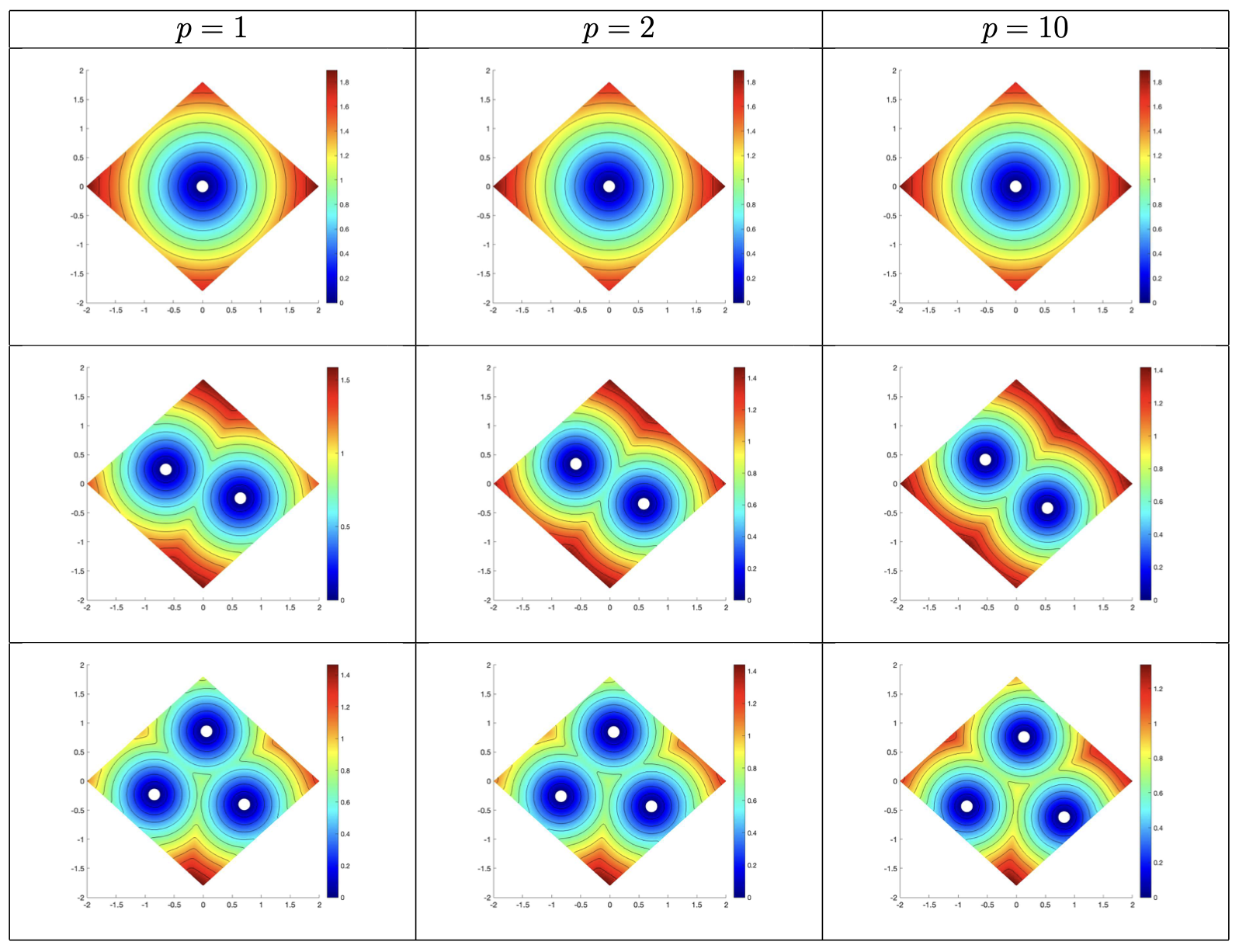}
    \caption{Optimal placement of $N\in \{1,2,3\}$ sensors in rhombus for different values of $p$.}
    \label{fig:sensors_4}
\end{figure}

In Figures \ref{fig:sensors_algeria}, \ref{fig:sensors_3} and \ref{fig:sensors_4}, we present some examples of the obtained optimal placement of $N\in \{1,2,3\}$ sensors for different values of $p$.


\
\newpage

\
\newpage

\ 
\newpage
\ 
\begin{remark}
Figure~\ref{fig:sensors_4} illustrates an intriguing symmetry-breaking phenomenon. Although one might naturally expect the optimal placement of two sensors to preserve the inherent symmetries of the domain -- namely, those of the rhombus -- the most efficient configuration in fact breaks this symmetry.
\end{remark}

\begin{remark}
 {   The optimal placement of the sensors is not necessarily unique. For example, if we consider the domain $\Omega$ to be a disk and take $N\ge 2$ sensors, then any rotation of an optimal placement with respect to the center of (the disk) $\Omega$ is also optimal. }
\end{remark}

\begin{remark}\label{rk:local}
 { As explained above, the use of gradient descent can lead to critical points in our case since we only use a first order sensitivity analysis. For example, if we consider two sensors, $\Omega$ to be a square and $p=1$, we notice that the method will converge to different critical points (with different values of the mean distances) depending on the choice of the initial placement of the sensors as shown in Figure \ref{fig:local} below:}
    \begin{figure}[ht]
    \centering
\begin{tabular}{|c|c|}
\hline
    \includegraphics[scale=.14]{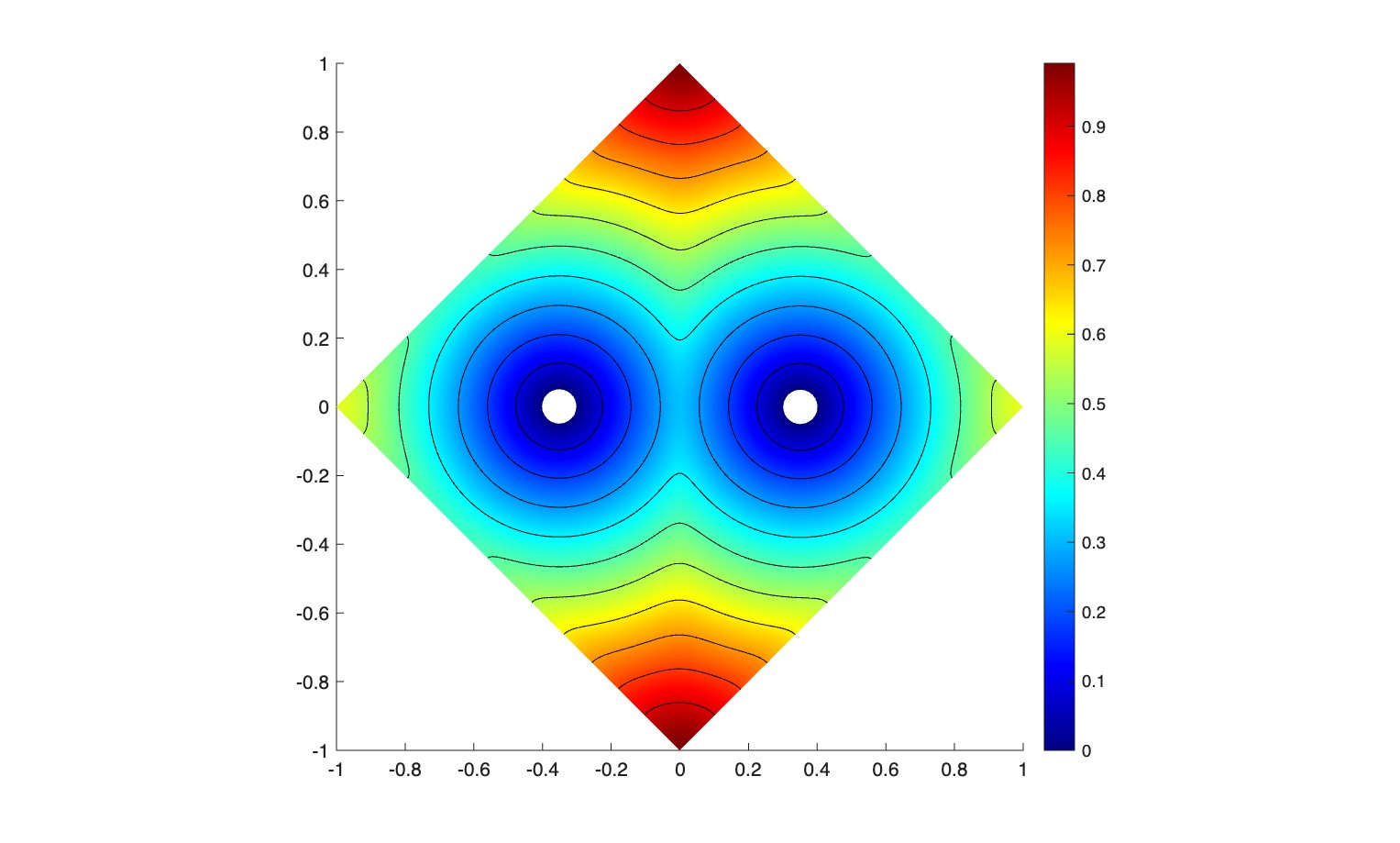}&
    \includegraphics[scale=.14]{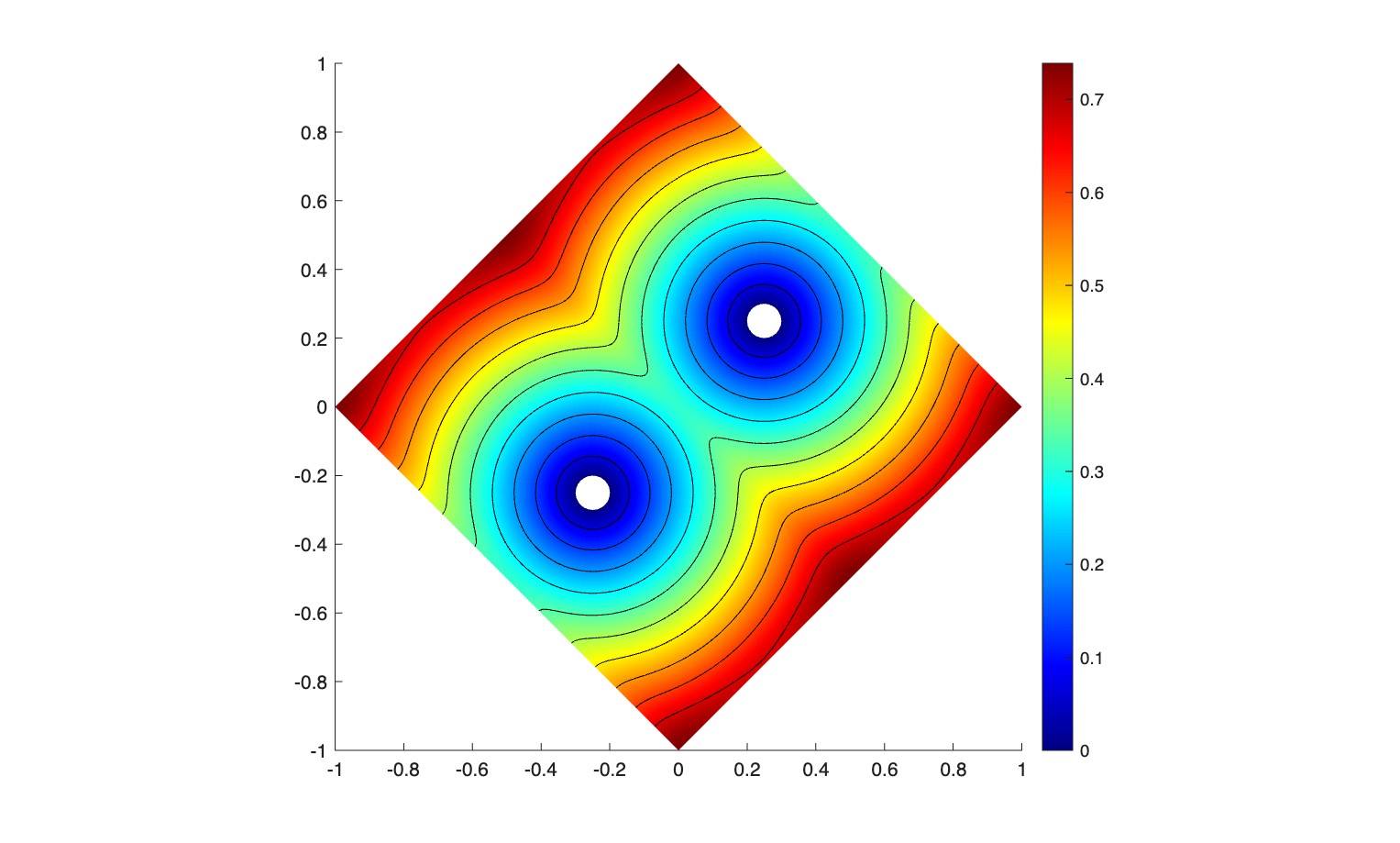}\\
\hline
    Mean distance $\approx 0.8066$&
    Mean distance $\approx 0.7947$\\
\hline
\end{tabular}
    \caption{Two critical configurations of the sensors obtained for $p=1$ by considering different initial positions: $(0.3,0)$ and $(-0.3,0)$ for the case on the left and $(0.1,0.1)$ and $(-0.3,-0.1)$ for the case on the right, which seems to be a global maximum.}
    \label{fig:local}
\end{figure}
\end{remark}

\section{Conclusion and perspectives}\label{s:conclusion}
The problem of optimal sensor placement has been addressed in a purely geometric setting, independent of the physical process under consideration and in the absence of PDE restrictions. The use of Varadhan's approximation theorem naturally leads to optimization problems constrained by the Laplacian. This allows us to apply the classical analytical and computational tools in optimal shape design with PDE criteria. 

\vspace{4mm}

\textbf{Conflict of Interest and data availability statements.} 
The authors declare that there is no conflict of interest. The datasets generated and analyzed during the current study are available from the corresponding author on reasonable request. \vspace{4mm}

\textbf{Acknowledgement.} The work of I. Ftouhi was supported by the Alexander von Humboldt Foundation via a Postdoctoral fellowship and also partially supported by the French National Research Agency under the France 2030 program, grant ANR-23-EXES-0005 (project Gardener/MASRA), and by the Occitanie Region.

The work of E. Zuazua was funded by the Alexander von Humboldt-Professorship program, the   ERC Advanced Grant CoDeFeL,  the Grants PID2020-112617GB-C22 KiLearn and TED2021-131390B-I00-DasEl of MINECO and PID2023-146872OB-I00-DyCMaMod of MICIU (Spain),  the European Union's Horizon Europe MSCA project ModConFlex, the Transregio 154 Project ``Mathematical Modelling, Simulation and Optimization Using the Example of Gas Networks'' of the DFG, the AFOSR 24IOE027 project, and the Madrid Government - UAM Agreement for the Excellence of the University Research Staff in the context of the V PRICIT (Regional Programme of Research and Technological Innovation).

\bibliographystyle{plain}
\bibliography{biblio}

\end{document}